\input amstex.tex
\input amsppt.sty   
\magnification 1200
\vsize = 8.5 true in
\hsize=6.2 true in
\nologo
\NoRunningHeads        
\parskip=\medskipamount
        \lineskip=2pt\baselineskip=18pt\lineskiplimit=0pt
       
        \TagsOnRight
        \NoBlackBoxes

        \topmatter
        Analysis of PDE
        \title
       Quasi-periodic solutions for nonlinear wave equations\\
       \\  
       {\it Solutions quasi-periodiques pour l'\'equation des ondes nonlin\'eaire}
         \endtitle

\author
         W.-M.~Wang        \endauthor        
\address
{CNRS and Department of Mathematics, Universit\'e Cergy-Pontoise, 95302 Cergy-Pontoise Cedex, France}
\endaddress
 \email
{wei-min.wang\@math.cnrs.fr}
\endemail
\abstract
We construct time quasi-periodic solutions to nonlinear wave equations on the torus in arbitrary dimensions.
All previously known results (in the case of zero or a multiplicative potential) seem to be limited to the circle. This generalizes the
method developed in the limit-elliptic setting in \cite{11} to the hyperbolic setting. The additional ingredient is a Diophantine property of algebraic numbers.

\noindent {\it R\'esum\'e. On construit des solutions quasi-p\'eriodiques en temps pour l'\'equation des ondes nonlin\'eaire sur le tore en dimension quelconque. 
Tous les r\'esultats pr\'ec\'edents se limitent au cercle. Cet article \'etend la m\'ethode d\'evelopp\'ee pour le cas limite-elliptique dans \cite{11} au cas hyperbolique. Le nouvel 
ingr\'edient est une propri\'et\'e diophantienne des nombres alg\'ebriques.}
\endabstract
     
        \bigskip\bigskip
        \bigskip
        \endtopmatter
          \bigskip
\document
\head{\bf 1. Introduction and statement of the Theorem}\endhead
We consider {\it real} valued solutions to the nonlinear wave  equation (NLW) on the $d$-torus $\Bbb T^d=[0, 2\pi)^d$:
$$
\frac{\partial^2v}{\partial t^2}-\Delta v+v+v^{p+1}+H(x,v)=0,\tag 1
$$
where $p\in\Bbb N$ and $p\geq 1$; considered as a function on $\Bbb R^d$, $v$ satisfies : $v(\cdot, x)=v(\cdot, x+2j\pi)$, $x\in [0, 2\pi)^d$ for all $j\in\Bbb Z^d$;
$H(x,v)$ is analytic in $x$ and $v$ and has the expansion:
$$H(x,v)=\sum_{m=p+2}^\infty \alpha_m(x) v^{m},$$
where $\alpha_m$ as a function on $\Bbb R^d$ is $(2\pi)^d$ periodic and real and analytic in a strip of width $\Cal O(1)$
for all $m$. The integer $p$ in (1) is {\it arbitrary}. 

Equation (1) can be rewritten as a first order equation in $t$. Let 
$$D=\sqrt {-\Delta+1}\tag 2$$ and 
$$u=(v, D^{-1}\frac{\partial v}{\partial t})\in\Bbb R^2.$$
Identifying $\Bbb R^2$ with $\Bbb C$, one obtains the corresponding first order equation 
$$i\frac{\partial u}{\partial t}=Du+D^{-1}[(\frac{u+\bar u}{2})^{p+1}+H(x, \frac{u+\bar u}{2})].$$

Using Fourier series, the solutions to the linear equation: 
$$i\frac{\partial u}{\partial t}=Du\tag 3$$
are linear combinations of eigenfunction solutions of the form:
$$e^{-i(\sqrt{j^2+1})t}e^{ij\cdot x},\quad j\in\Bbb Z^d,$$
where $j^2=|j|^2$ and $\cdot$ is the usual inner product. These solutions are either periodic or quasi-periodic in time. 

The main purpose of this note is to announce the following new result, namely, under appropriate genericity conditions on the Fourier frequencies, to be detailed  in sect. 2, 
a class of quasi-periodic solutions to the linear wave equation (3) bifurcates to quasi-periodic solutions of the NLW in (1). 
We note that when $p\geq\frac{4}{d-2}$ ($d\geq 3$, $H=0$), global solutions to (1) do not seem to 
be known in general. 

Under the assumption that $H$ is a polynomial in $u$, $\bar u$, $e^{ix_k}$ and $e^{-ix_k}$, $k=1, 2, ..., b$, $x_k\in[0, 2\pi)$, below is the precise statement. 

\proclaim
{Theorem} 
Assume that $$v^{(0)}(t, x)={\text Re }\sum_{k=1}^b a_k e^{-i(\sqrt{j_k^{2}+1})t}e^{ij_k\cdot x}$$ is generic, satisfying the genericity conditions (i-iii), $a=\{a_k\}\in (0,\delta]^b=\Cal B(0, \delta)$ and $p$ even. 
Assume that $b>C_p d$. Then for all $0<\epsilon<1$,  there exists $\delta_0>0$, such that 
for all $\delta\in (0, \delta_0)$,  there is a Cantor set $\Cal G\subset \Cal B(0, \delta)$ with 
$$\text{meas }\Cal G/\delta^b\geq 1-\epsilon.$$ 
For all $a\in\Cal G$, there is a quasi-periodic solution of $b$ frequencies to the nonlinear wave equation (1):
$$v(t, x)={\text Re}\sum a_k e^{-i{\omega_k}t}e^{ij_k\cdot x}+o(\delta^{3/2}),$$
with basic frequencies $\omega=\omega(a)=\{\omega_k(a)\}^b_{k=1}$ satisfying 
$$\omega_k=\sqrt{j_k^2+1}+\Cal O(\delta^{p}),$$
and the amplitude-frequency map $a\mapsto\omega (a)$ is a diffeomorphism. 
The remainder $o(\delta^{3/2})$ is in a Gevrey norm on $\Bbb T^{b+d}$.
\endproclaim

\noindent{\it Remark.}  One views quasi-periodic solutions of $b$ frequencies as periodic solutions on a $b$-dimensional torus in ``time". Hence the use of Gevrey norms on $\Bbb T^{b+d}$.
The condition of large $b$, namely $b>C_pd$, is imposed in order that certain determinants are not identically zero, as in \cite{11}. It cannot be excluded that this condition could be improved after 
more  technical work. Contrary to \cite{11}, however, aside from the genericity conditions, this is the only other condition needed to prove the Theorem.  This is because the genericity condition (ii) below 
dictates that $\omega$ is Diophantine, see (4) below; moreover the mass term $1$ in the wave operator $D$, introduces curvature, cf. \cite{12}. The polynomial restriction on $H$ is technical, the result most likely remains valid for analytic $H$.
 \smallskip 
This Theorem appears to be the first general existence results on quasi-periodic solutions to the NLW  in (1) in arbitrary dimensions. Previously quasi-periodic solutions only seem to have been constructed in one dimension with positive mass $m$. In that case, the linear wave equation:
$$
\frac{\partial^2v}{\partial t^2}-\frac{\partial^2v}{\partial x^2}+mv=0,
$$ 
gives rise to an eigenvalue set $\{\sqrt {j^2+m}, j\in\Bbb Z \}$ close to the set of integers, see \cite{4, 9} and \cite{6, 8, 13}
in a related context. For almost all $m$, this 
set is linearly independent over the integers. This property does not have higher dimensional analogues and seems to have been a serious obstacle.
(The time periodic case is special and solutions have been constructed in higher dimensions in \cite{3}.) 

To insert into a more general context, it is known, cf. \cite {7},  for example,  that the spectrum of an elliptic first-order operator on generic 
compact manifold is {\it dense}. Since the spectrum of the wave operator $D$ in (2) on the flat torus is the set 
$$\{\sqrt{j^2+1}|\, j\in\Bbb Z^d\},$$
which becomes dense in $d\geq 2$, it can serve as a model example of the generic case.  
\bigskip
\head{\bf 2. The generic linear solutions }\endhead  
As in \cite{11}, the Fourier support: 
$$\tilde\Gamma=\bigcup_{r=1}^{R} \text{supp }[(u^{(0)}+{\bar u}^{(0)})^{*pr}]\backslash\{(0, 0)\}\subset\Bbb Z^{b+d},$$
for some $R>0$, plays an essential role. Using  $(\nu, \eta)$ to denote an element of $\tilde\Gamma$: $(\nu, \eta)\in \tilde\Gamma$, then as in sect. 2 of \cite{11}, there is the following relation:
$$\eta=-\sum_{k=1}^b\nu_kj_k,$$
and consequently $\eta=\eta(j_1, j_2,..., j_b)$ is viewed as a function from $(\Bbb Z^d)^b$ to $\Bbb Z^d$. 

For the purposes below, $(\nu, \eta)$ is, however,  viewed as a point in $\Bbb Z^{b+d}$.

\noindent{\it Definition.} $u^{(0)}$ of $b$ frequencies $j_1, j_2, ..., j_b\in \Bbb Z^d$ is {\it generic} if the following three conditions are satisfied:
\item{(i)} Let $\bar d=\min (d, b)$. Any $\bar d$ vectors in the set $\{j_k \}_{k=1}^b$ are linearly independent. 
For all $j_k$, $k=1, 2, ..., b$, define the set of differences 
$$J_k=\{j_{k'}-j_k|k'=1, ..., b, k'\neq k\}.$$ 
If $b\geq d+1$, any $d$ vectors in $J_k$ are linearly independent.
(If $b\leq d$, there is no condition (i).)
\item{(ii)} The integers $(j^2_k+1)$, $k=1$, $2$, ..., $b$, are distinct:
$$1<j_1^2+1<j_2^2+1<\cdots<j_b^2+1,$$
and square free.
\item{(iii)} For all $k=1, 2, ..., b$, $m\in[-p, p]\backslash\{0\}$, consider the set of $\eta$ of the form 
$$\aligned \eta=&mj_k-m_\ell j_\ell\\
:=&\eta(m_\ell)\neq 0,\endaligned$$
where $\ell=1, 2, ..., b$, $m_\ell \in\Bbb Z$, $|m_\ell|\leq 2p(d+1)$. For each $\eta$, define $L$ to be 
$$L=2m\eta\cdot j_k+(m^2-m_\ell^2):=L(m_\ell).$$
Denote by $P(m_\ell)$ the corrsponding $d$-dimensional hyperplane in $\Bbb R^d$:
$$2\eta\cdot j+L=0,$$ 
where $\eta=\eta(m_\ell)$ and $L=L(m_\ell)$. 

Let $\sigma$ be any set of $(d+1)$ $m_\ell$, $\ell=1, 2, ..., b$,  $m_\ell\in\Bbb Z$, $|m_\ell|\leq 2p(d+1)$, such that there exists $m_\ell\in\sigma$, $m_\ell\neq \pm m$, then 
$$\bigcap_{\sigma} P(m_\ell)=\emptyset.$$
\smallskip
\noindent{\it Remark.} In lieu of condition (ii), one may take $(j_k^2+1)$ to be multiples of {\it distinct} 
square free integers -- the proof of the Theorem is the same.
\smallskip
The following indicates that the above three conditions are viable.
\proclaim {Lemma} There are infinite numbers of $j_1$, $j_2$, ..., $j_b\in\Bbb Z^d$ which satisfy the genericity conditions (i-iii).
\endproclaim
\noindent The proof of (i, iii) for the Lemma is similar to the corresponding one in \cite{11}. Clearly there are infinite numbers of integers satisfying (ii).  In fact given $N\in\Bbb N$,
the number of square free integers less than or equal to $N$ is asymptotically, 
$$= \frac {6}{\pi^2}N+\Cal O(\sqrt N).$$
\hfill $\square$
\smallskip 
It follows from basic algebra that condition (ii) implies the usual (linear) Diophantine property: 
$$\Vert \sum_{k=1}^bn_k\omega_k\Vert_{\Bbb T}\neq 0,$$
where $n_k\in\Bbb Z$, $\sum_k|n_k|\neq 0$;
as well as the quadratic Diophantine property: 
$$\Vert \sum_{k, \ell; k< \ell}n_{k\ell}\omega_k\omega_\ell \Vert_{\Bbb T}\neq 0,$$
where $n_{k \ell}\in\Bbb Z$, $\sum |n_{k\ell}|\neq 0$.
So $$\Vert \sum_{k=1}^bn_k\omega_k\Vert_{\Bbb T}\geq c |n|^{-\alpha}\tag 4$$
for some $\alpha >0$,  and $|n|=\sum |n_{k}|\neq 0$ and 
$$\Vert \sum_{k, \ell, k<\ell} n_{k \ell} \omega_k\omega_{\ell}\Vert _{\Bbb T} \geq c' |n|^{-\beta},\tag 5$$
for some $\beta >0$,  and $|n|=\sum |n_{k \ell}|\neq 0$, cf. \cite{10}. 

Consequently, 
$$\Vert \sum_{k, \ell} n_{k \ell} \omega_k\omega_{\ell}\Vert _{\Bbb T}=0\, \Longleftrightarrow \sum_{k, \ell} n_{k \ell} \omega_k\omega_{\ell}=n_{kk} \omega^2_k,\tag 6$$
for some $k\in\{1,2, ..., b\}$.
The expression above of linear combinations of products of pairs of eigenvalues appears as the principal symbol of an 
appropriate linearized operator, which is the ``divisor" in the problem. The above ``zero-divisors" 
are dealt with as in \cite{11}, which essentially uses the sub-principal symbol  of the linearized 
operator to control the small eigenvalues under the genericity conditions (i, iii). Combining with the small-divisor estimates in (5), one is then able to achieve 
amplitude-frequency modulation. One can then adapt the analysis construction of Bourgain in Chap. 20 of \cite{5} to complete the proof of the Theorem as in \cite{11}.
\hfill$\square$

\noindent {\it Remark.} In Chap. 20 of \cite{5}, when $\omega$ is a Fourier multiplier - an external parameter, nonlinear polynomial conditions on $\omega$
are imposed in oder to achieve separation properties. It can be shown that, in fact, quadratic polynomial conditions suffice by using antisymmetry property of the 
determinants, cf. \cite {1, 2, 12}.  When $\omega$ is {\it fixed}, in view of (6), the corresponding condition simply {\it cannot} be imposed -- instead the genericity conditions (i-iii) are 
used to bypass this obstruction to invertibility.

\enddocument
\end